\documentclass[reqno,11pt]{article}
\usepackage{amssymb}
\usepackage{amsmath}
\usepackage{amsthm}
\usepackage{geometry} 

\newtheorem{Theorem}{Theorem } [section]
\newtheorem{lemma}[Theorem]{Lemma}

\numberwithin{equation}{section}

\newcommand{\vertiii}[1]{{\vert\kern-0.25ex\vert\kern-0.25ex\vert #1 
    \vert\kern-0.25ex\vert\kern-0.25ex\vert}}

\DeclareMathOperator{\supp}{supp}

\date{}  

\title{A note on the uniqueness properties of solutions for the Schrödinger-Korteweg de Vries system}
\author{Eddye Bustamante, Jos\'e Jim\'enez Urrea and Jorge Mej\'{\i}a}

\author{
    Eddye Bustamante \\
    \textit{Universidad Nacional de Colombia, Sede Medell\'in} \\
    \texttt{eabusta0@unal.edu.co} 
    \and
    José Jiménez Urrea \\
    \textit{Universidad Nacional de Colombia, Sede Medell\'in} \\
    \texttt{jmjimene@unal.edu.co} 
    \and
    Jorge Mej\'ia \\
    \textit{Universidad Nacional de Colombia, Sede Medell\'in} \\
    \texttt{jemejia@unal.edu.co} 
}

\begin{document}

%
%

\maketitle

\begin{abstract} In this work we prove that if $(u_i,v_i)$, $i=1,2$, are smooth enough solutions of the coupled Schrödinger-Korteweg-de Vries system
\begin{align*}
\left. \begin{array}{rl}
i u_t+\partial_x^2 u &\hspace{-2mm}=\beta uv - |u|^2 u,\\
\partial_t v + \partial_x^3 v &\hspace{-2mm}=\gamma \partial_x |u|^2-\frac12\partial_x (v^2)
\end{array} \right\}
\end{align*}
with appropriate decay at infinity such that at two different times $t_0=0$ and $t_1=1$ satisfy that
$$u_1(0)-u_2(0),u_1(1)-u_2(1),v_1(0)-v_2(0),v_1(1)-v_2(1)\in H^1(e^{ax^{2}}dx),$$
for $a>0$ big enough, then $u_1=u_2$ and $v_1=v_2$.\\

(Let us recall that $f\in H^1(e^{ax^{2}} dx)$ iff $f\in L^2(e^{ax^{2}}dx)$ and $\partial_x f\in L^2(e^{ax^{2}}dx)$).
\end{abstract}

\setcounter{page}{1}
\pagenumbering{arabic}

\section{Introduction}

\footnotetext[1]{\textbf{2020 Mathematics Subject Classification:} 35Q53}

In this article, we consider the initial value problem (IVP) associated with the Schrödinger-Korteweg-de Vries system (NLS-KdV) on the line
\begin{align}
\left. \begin{array}{rl}
i u_t+\partial_x^2 u &\hspace{-2mm}=\beta uv - |u|^2 u, \quad x\in\mathbb R,\quad t\in\mathbb R,\\
\partial_t v + \partial_x^3 v &\hspace{-2mm}=\gamma \partial_x |u|^2-\frac12\partial_x (v^2), \quad x\in\mathbb R,\quad t\in\mathbb R,\\
u(x,0)=u_0(x), & v(x,0)=v_0(x),
\end{array} \right\}\label{maineq}
\end{align}

where $u=u(x,t)$ is a complex valued function, $v=v(x,t)$ is a real valued function and $\beta,\gamma$ are real constants.\\ 

The system \eqref{maineq} has been widely studied in various mathematical and physical contexts. From the physics point of view, the NLS-KdV system arises in the study of the interaction between short and long waves, as well as in the analysis of resonance phenomena in capillary-gravity water waves (see \cite{FO1983},  \cite{HIMN1974} and \cite{KKS1975}).\\

As shown in \cite{BB1983}, the system \eqref{maineq} is not completely integrable. Its well-posedness, however, has been extensively studied in the literature. Below, we mention some important results for the case where the initial data $(u_0,v_0)\in H^k(\mathbb R)\times H^s(\mathbb R)$.\\

In \cite{T1993}, Tsutsumi proved that the system \eqref{maineq} is globally well-posed in $H^{s+1/2}(\mathbb R)\times H^s(\mathbb R)$ for $s\in \mathbb Z^+$. Later, Bekiranov, Ogawa and Ponce, in \cite{BOP1997}, extended Tsutsumi result to real values $s>0$, using the restriction method of the Fourier transform (see \cite{B1993}). Generalizing the admisible indices for the choice of the initial data, Corcho and Linares, in \cite{CL2007}, proved the local well-posedness of the system \eqref{maineq} for  initial data $(u_0,v_0)\in H^k(\mathbb R)\times H^s(\mathbb R)$ with $k\geq 0$, $s>-\frac34$, and
$$
\left\{
\begin{array}{lr}
k-1\leq s\leq 2k-\frac12& \text{ if } 0\leq k\leq \frac12,\\
k-1\leq s\leq k+\frac12 & \text{ if } k>\frac12.
\end{array}
\right.
$$

Improving the multilinear estimates obtained in \cite{CL2007}, Wu, in \cite{W2010}, extended the previous result to the region
$$A_0:=\{(k,s)\in\mathbb R^2: k\geq 0, s>-\frac34,s<4k,-1<k-s<2\}.$$

More recently, Correia, Linares and Drumond Silva established in \cite{CLDS2025} an optimal result by considering the region
$$A=\{(k,s)\in\mathbb R^2: k\geq 0, s>-\frac34, s<4k, -2<k-s<3\},$$
proving local well-posedness in $H^k(\mathbb R)\times H^s(\mathbb R)$ for $(k,s)\in A$, and ill-posedness when $(k,s)\notin \overline A$. This result was obtained through a novel approach, not by modifying the funcional spaces to apply fixed point methods, but instead by introducing the concept of integrated-by-parts strong solutions, which can be viewed as a generalization of the classical strong solution concept in the sense of Duhamel's formula.\\

The aim of this work is to establish unique continuation properties for the system \eqref{maineq}. In general terms, this refers to the fact that if the difference of two sufficiently regular solutions of \eqref{maineq} decays as $e^{-ax^{2}}$ at two different times, then the solutions must coincide. Results in this nature have been obtained for each of the equations in the system \eqref{maineq}. In \cite{EKPV2006}, Escauriaza, Kenig, Ponce and Vega studied unique continuation properties of solutions for Scrödinger-type equations of the form
\begin{align}
i\partial_t u +\Delta u=Vu,\quad (x,t)\in \mathbb R^n\times \mathbb R,\label{1.2}
\end{align}

and established sufficient conditions on the behavior of the solution $u$ at two different times $t_0=0$ and $t_1=1$, which guarantee that $u\equiv 0$ is the unique solution of \eqref{1.2}. Their result is the following.\\

\textbf{Theorem A.} \textit{Let $u\in C([0,1];H^2(\mathbb R^n))$ be a strong solution of the equation \eqref{1.2} in $\mathbb R^n\times[0,1]$ with $V:\mathbb R^n\times[0,1]\to\mathbb C$ such that $V\in L^\infty(\mathbb R^n\times[0,1])$ and $\nabla_x V\in L^1_t([0,1];L^\infty(\mathbb R^n_x))$. Then there exists $a_0>0$ such that for $a\geq a_0$, if $u(0),u(1)\in H^1(e^{a|x|^2}dx)$, and
\begin{align}
\lim_{r\to\infty} \|V\|_{L^1_t L^\infty_x(\{|x|>r\})}=\lim_{r\to\infty}\int_0^1 \sup_{|x|>r} |V(x,t)|dt=0,\label{1.4}
\end{align}
then $u\equiv 0$.}\\

\newpage
In \cite{EKPV2007}, Escauriaza, Kenig, Ponce and Vega, study uniqueness properties of solutions of the k-generalized Korteweg-de Vries equations of the form
\begin{align}
\partial_t v +\partial_x^3 v=-v^k \partial_x v,\label{1.5}
\end{align}

and establish sufficient conditions on the behavior of the difference $v_1-v_2$ of two solutions $v_1,v_2$ of \eqref{1.5} at two different times $t_0=1$ and $t_1=1$, which guarantee that $v_1=v_2$. Their result can be formulated as follows.\\

\textbf{Theorem B.} \textit{Let $v_1,v_2\in C([0,1];H^3(\mathbb R)\cap L^2(|x|^2 dx))$ be strong solutions of \eqref{1.5}. If $v_1(0)-v_2(0)$ and $v_1(1)-v_2(1)\in H^1(e^{a|x|^{3/2}} dx)$ for any $a>0$, then $v_1\equiv v_2$.}\\

The above theorems can be used directly to prove a uniqueness property of the system \eqref{maineq} as follows:\\

Let $(u,v)\in C([0,1]; H^{7/2}(\mathbb R_x)\cap L^2(|x|^7 dx))\times C([0,1]; H^3(\mathbb R_x)\cap L^2(|x|^3 dx))$ be a solution of the system \eqref{maineq}. Then there exists $a_0>0$ such that for some $a\geq a_0$ if $u(0),u(1)\in H^1(e^{ax^2}dx)$, and $v(0),v(1)\in H^1(e^{a|x|^{3/2}} dx)$, then $u\equiv 0$ and $v\equiv 0$.\\

This result follows from Theorems A and B, because $u$ satisfies the equation
\begin{align*}
i\partial_t u + \partial_x^2 u=Vu,
\end{align*}
where $V:=(\beta v-|u|^2)\in L^\infty(\mathbb R\times[0,1])$ is such that $\partial_x V \in L^1_t([0,1];L^\infty_x(\mathbb R))$ and $V$ satisfies \eqref{1.4}. In consequence, by Theorem A, $u\equiv 0$. Therefore, $v$ satisfies the equation
\begin{align*}
\partial_t v + \partial_x^3 v=-\frac12 \partial_x(v^2),
\end{align*}

and by Theorem B, $v\equiv 0$.\\

However, if $(u_1,v_1), (u_2,v_2)\in C([0,1];H^{7/2}(\mathbb R_x)\cap L^2(|x|^7 dx))\times C([0,1];H^3(\mathbb R_x)\cap L^2(|x|^3 dx))$ are two solutions of the system \eqref{maineq} such that
\begin{align}
u_1(0)-u_2(0),u_1(1)-u_2(1)\in H^1(e^{ax^2} dx),\text{ and }v_1(0)-v_2(0),v_1(1)-v_2(1)\in H^1(e^{a|x|^{3/2}}dx)\label{1.6}
\end{align}
for $a\geq a_0$, being $a_0$ certain positive number, we cannot apply the same reasoning above to the difference $u:=u_1-u_1$, $v:=v_1-v_2$, because in this case the equation satisfied by $u$ is
\begin{align*}
i\partial_t u + \partial_x^2 u = (\beta v - |u_1|^2 - u_2 \overline u_1) u + (\beta u_2 v - u_2^2 \overline u),
\end{align*}

which is not of the form \eqref{1.2}, to apply Theorem A.\\

For this reason, in this article we have to develop a lower estimate (Theorem \ref{le_th} below) that allows us to prove a uniqueness principle for two solutions $(u_1,v_1)$ and $(u_2,v_2)$ of the system \eqref{maineq}, but with more restrictive hypotheses about $u:=u_1-u_2$ and $v:=v_1-v_2$, at two different times $t_0=0$ and $t_1=1$.\\

The proof of the uniqueness principle for the system \eqref{maineq} follows the same ideas of \cite{EKPV2007} based upon the confrontation of Carleman estimates and the lower estimate. In this case the Carleman estimates, we use, are those well-known Carleman estimates for the linear Schrödinger and KdV equations (see Section \ref{SectionCarleman}).\\

Our results in this note can be stated, precisely, as follows:
\begin{Theorem}\label{le_th} (Lower estimate). Let $(u_i,v_i)\in C([0,1];H^{7/2}(\mathbb R)\cap L^2(|x|^7 dx))\times C([0,1];H^3(\mathbb R)\cap L^2(|x|^3 dx))$ $i=1,2$, be solutions of the Schrödinger-KdV system \eqref{maineq} and define $\tilde u:=u_1-u_2$, $\tilde v:=v_1-v_2$. For $r\in(0,\frac12)$ let $Q_r:=\{(x,t): 0\leq x\leq 1,t\in[r,1-r]\}$ and suppose that $\|\tilde u\|_{L^2(Q_{r_0})}+\|\tilde v\|_{L^2(Q_{r_0})}>0$ for some $r_0\in(0,\frac12)$.\\

For $R>0$, let us define
\begin{align}
\notag A_R(\tilde u,\tilde v)&:=\left( \int_0^1 \int_{R-1}^R (|\tilde u|^2+|\partial_x \tilde u|^2) dx dt \right)^{1/2} + \left( \int_0^1 \int_{R-1}^R (|\tilde v|^2+|\partial_x \tilde v|^2 + |\partial_x^2 \tilde v|) dx dt \right)^{1/2}\\
&\equiv A_{s,R}(\tilde u) + A_{k,R}(\tilde v). \label{3.4}
\end{align}
Then there are constants $C>0$, $R_0\geq 2$ and $r_1\in (0,r_0)$ such that for all $r\in(0,r_1)$ there is $\overline C=\overline C(r)>0$ for which
\begin{align}
\|\tilde u\|_{L^2(Q_r)} + \|\tilde v\|_{L^2(Q_r)} \leq C e^{9 \overline C R^{2}} A_R(\tilde u,\tilde v),\quad \text{for all }R\geq R_0.\label{3.5}
\end{align}

\end{Theorem}

The main result concerning uniqueness properties of the Schrödinger-Korteweg-de Vries system is the following.

\begin{Theorem}\label{main_th} Suppose that $(u_i,v_i)\in C([0,1];H^{7/2}(\mathbb R)\cap L^2(|x|^7 dx))\times C([0,1]; H^3(\mathbb R)\cap L^2(|x|^3 dx))$, $i=1,2$, are solutions of the system \eqref{maineq}. Then there exists a universal constant $a_0>0$, such that if for some $a\geq a_0$
\begin{align}
u_1(0)-u_2(0),\, u_1(1)-u_2(1),\, v_1(0)-v_2(0),\, v_1(1)-v_2(1)\in H^1(e^{ax^{2}} dx),\label{4.2}
\end{align}
then $u_1=u_2$ and $v_1=v_2$.\\
(Let us recall that $f\in H^1(e^{a x^{2}} dx)$ if and only if $f\in L^2(e^{a x ^{2}}dx)$ and $\partial_x f\in L^2(e^{a x^{2}} dx)$ ).
\end{Theorem}

The existence of solutions for the Schrödinger-Korteweg-de Vries system \eqref{maineq} with the regularity and decay required in Theorem \ref{main_th} is guaranteed thanks to a persistence result obtained by Linares and Palacios in \cite{LP2019} (see Theorem 1.2 of that article).\\

The question that remains open is the following: can the uniqueness principle be established under the hypothesis \eqref{1.6} instead of \eqref{4.2}?\\

This paper is organized as follows: in Section \ref{SectionCarleman} we state the Carleman type estimates for the KdV and the Scrödinger equations that will be used in the proof of Theorem \ref{main_th}. Section \ref{SectionLE} is devoted to prove the lower estimate stated in Theorem \ref{le_th}. Finally, in Section \ref{SectionMR}, we proof Theorem \ref{main_th}.\\

Throughout the paper the letter $C$ will denote diverse constants, which may change from line to line, and whose dependence on certain parameters is clearly established in all cases.

\section{Carleman type estimates}\label{SectionCarleman}

In this section, we present the Carleman type estimates that will be used in the proof of Theorem \ref{main_th}. These estimates are inspired by the works of Kenig, Ponce and Vega  \cite{KPV2003}, as well as Escauriaza, Kenig, Ponce and Vega \cite{EKPV2006} and \cite{EKPV2007}. However, our results are established in $L^pL^q$ spaces with simpler indices $p$ and $q$, by means of elementary Fourier transform techniques, following the approach used, for instance, by Bustamante, Isaza and Mejía in \cite{BIM2012}.\\

Concerning the linear KdV equation, we have the following result.
\begin{lemma}\label{carleman_kdv} Let $D:=\{(x,t): x\in\mathbb R, t\in[0,1]\}$, and let $w\in C([0,1];H^3(\mathbb R))\cap C^1([0,1];L^2(\mathbb R))$ be such that, for some $M>0$, $\supp w(t)\subset [-M,M]$ for every $t\in[0,1]$. Then, for every $\lambda>2$,
\begin{align}
\|e^{\lambda x}w\|_{L^\infty_t L^2_x(D)} \leq \|e^{\lambda x} w(0)\|_{L^2_x(\mathbb R)} + \|e^{\lambda x}w(1)\|_{L^2_x(\mathbb R)} + \|e^{\lambda x} (w_t+\partial_x^3 w)\|_{L^1_t L^2_x(D)}. \label{carl_kdv_1}
\end{align}
Moreover, there exists $C>0$, independent of $w$ and $M$, such that, for every $\lambda>2$,
\begin{align}
\|e^{\lambda x} \partial_x w\|_{L^\infty_x L^2_t(D)} & \leq C \lambda^2 ( \|e^{\lambda x} w(0)\|_{L^2_x(\mathbb R)} + \|e^{\lambda x}w(1)\|_{L^2_x(\mathbb R)} ) + C \|e^{\lambda x} (w_t+\partial_x^3 w)\|_{L^1_x L^2_t(D)}, \label{carl_kdv_2}\\
\|e^{\lambda x} \partial^2_x w\|_{L^\infty_x L^2_t(D)} & \leq C \lambda^2 ( \| J(e^{\lambda x} w(0))\|_{L^2_x(\mathbb R)} + \|J(e^{\lambda x}w(1))\|_{L^2_x(\mathbb R)} ) + C \|e^{\lambda x} (w_t+\partial_x^3 w)\|_{L^1_x L^2_t(D)}, \label{carl_kdv_3}
\end{align}
where $\widehat{Jg}(\xi)=(1+\xi^2)^{1/2} \widehat g(\xi)$.
\end{lemma}

Now we turn to the linear Schrödinger equation, for which the corresponding result is the following.
\begin{lemma}\label{carleman_schr} Let $D:=\{(x,t): x\in\mathbb R, t\in[0,1]\}$, and let $w\in C([0,1];H^2(\mathbb R))\cap C^1([0,1];L^2(\mathbb R))$ be such that, for some $M>0$, $\supp w(t)\subset [-M,M]$ for every $t\in[0,1]$. Then, for every $\lambda>2$,
\begin{align}
\|e^{\lambda x}w\|_{L^\infty_t L^2_x(D)} \leq \|e^{\lambda x} w(0)\|_{L^2_x(\mathbb R)} + \|e^{\lambda x}w(1)\|_{L^2_x(\mathbb R)} + \|e^{\lambda x} (iw_t+\partial_x^2 w)\|_{L^1_t L^2_x(D)}. \label{carl_schr_1}
\end{align}
Moreover, there exists $C>0$, independent of $w$ and $M$, such that, for every $\lambda>2$,
\begin{align}
\|e^{\lambda x} \partial_x w\|_{L^\infty_x L^2_t(D)} & \leq C \lambda^2 ( \| J(e^{\lambda x} w(0))\|_{L^2_x(\mathbb R)} + \|J(e^{\lambda x}w(1))\|_{L^2_x(\mathbb R)} ) + C \|e^{\lambda x} (iw_t+\partial_x^2 w)\|_{L^1_x L^2_t(D)}. \label{carl_schr_3}
\end{align}
\end{lemma}

\section{Proof of Theorem \ref{le_th} (lower estimate)}\label{SectionLE} Following the ideas in \cite{I1993}, \cite{EKPV2006} and \cite{EKPV2007}, it can be proved the following two lemmas, which express boundedness properties of the inverses of the operators associated to the linear parts of the Schrödinger (Lemma \ref{le_schr}) and KdV (Lemma \ref{le_kdv}) equations.

\begin{lemma}\label{le_schr} Let $\phi:[0,1]\to \mathbb R$ be a $C^\infty$ function and let $D:=\{(x,t): x\in\mathbb R, t\in[0,1]\}$. Let us assume that $R>1$ and define
\begin{align}
\psi(x,t):=\alpha \left( \frac xR +\phi(t) \right)^2, \label{3.1}
\end{align}
being $\alpha$ a free parameter. Then, there is $\overline C:=\max\{\|\phi'\|_{L^\infty},\|\phi'\|_{L^\infty}^2,\|\phi''\|_{L^\infty},1\}>0$ such that the inequality
\begin{align}
\frac{\alpha^{3/2}}{R^2} \|e^{\psi}g\|_{L^2(D)}+\frac{\alpha^{1/2}}{R} \|e^{\psi} \partial_x g\|_{L^2(D)}\leq \sqrt2 \|e^{\psi}(i\partial_t+\partial_x^2) g\|_{L^2(D)} \label{3.2}
\end{align}
holds if $\alpha\geq \overline C R^2$ and $g\in C([0,1]; H^2(\mathbb R))\cap C^1([0,1]; L^2(\mathbb R))$
is a complex valued function such that
\begin{enumerate}
\item[(i)] $g(0)=g(1)=0$;
\item[(ii)] there is $M>0$ such that $\supp g(t)\subset[-M,M]$ for all $t\in[0,1]$.
\item[(iii)] $\supp g(\cdot_t)(\cdot_x)\subset \{(x,t):|\frac xR+\phi(t)|\geq 1\}$.
\end{enumerate}
\end{lemma}

\begin{lemma}\label{le_kdv} Let $\phi,\psi$, and $D$ as in Lemma \ref{le_schr}. There is $\overline C:=\max\{\|\phi'\|_{L^\infty}, \|\phi''\|_{L^\infty},1\}>0$ such that the inequality
\begin{align}
\frac{\alpha^{5/2}}{R^3} \|e^{\psi}g\|_{L^2(D)}+\frac{\alpha^{3/2}}{R^2} \|e^{\psi} \partial_x g\|_{L^2(D)}\leq \sqrt2 \|e^{\psi}(\partial_t+\partial_x^3) g\|_{L^2(D)} \label{3.3}
\end{align}
holds if $\alpha\geq \overline C R^{3/2}$ and $g\in C([0,1]; H^3(\mathbb R))\cap C^1([0,1]; L^2(\mathbb R))$
is a real valued function satisfying conditions (i), (ii) and (iii) in Lemma \ref{le_schr}.
\end{lemma}

\textit{Proof of Theorem \ref{le_th}}. Proceeding as in \cite{BIM2013}, let us take $\phi \in C^{\infty}([0,1])$ such that $\phi=0$ in $[0,\frac r2]\cup [1-\frac r2,1]$, $\phi=4$ in $[r,1-r]$, $\phi$ is increasing in $[\frac r2,r]$ and $\phi$ is decreasing in $[1-r,1-\frac r2]$, and $\|\phi'\|_{L^\infty}\geq \frac{1}{r}$ ($r\in(0,\frac12)$ to be determined later). Let $\mu\in C^{\infty}(\mathbb R)$ be an increasing function with $\mu=0$ in $(-\infty,2]$ and $\mu=1$ in $[3,\infty)$. Let also consider $\theta \in C^\infty(\mathbb R)$ such that $\theta=1$ if $x\leq R-1$, $\theta=0$ if $x\geq R$ and such that $|\theta'|,|\theta''|,|\theta'''|\leq C$, where $C$ is a constant independent of $R$.\\

We will apply Lemma \ref{le_schr} to the function
\begin{align}
\tilde g(t)(x):=\mu\left(\frac xR + \phi(t)\right) \theta(x) \tilde u(t)(x),\label{3.6}
\end{align}
and Lemma \ref{le_kdv} to the function
\begin{align}
g(t)(x):=\mu\left(\frac xR + \phi(t)\right) \theta(x) \tilde  v(t)(x).\label{3.7a}
\end{align}

For the sake of simplicity we will write $\tilde g=\mu\theta\tilde u$ and $g=\mu\theta\tilde v$.\\

From the definition of $g$ and the fact that $\tilde v$ satisfies the equation
\begin{align}
\partial_t \tilde v + \partial_x^3 \tilde v + (v_1+v_2) \partial_x \tilde v + \partial_x(v_1+v_2) \tilde v - \gamma (\overline u_1 \partial_x \tilde u + u_2\partial_x \overline{\tilde u}) - \gamma (\partial_x u_1 \tilde u + \partial_x u_2 \overline{\tilde u})=0, \label {3.7}
\end{align}

it can be seen that
\begin{align}
\partial_t g + \partial_x^3 g = -\mu\theta(v_1+v_2) \partial_x \tilde v - \mu\theta \partial_x(v_1+v_2)\tilde v + F_1+F_2+F_3,\label{3.8}
\end{align}

where
\begin{align}
F_1&:=\mu[\theta''' \tilde v+3\theta'' \partial_x \tilde v+3\theta' \partial_x^2 \tilde v],\label{3.9}\\
F_2&:=\frac{\mu'}R 3\theta \partial_x^2 \tilde v + \left[\frac3{R^2}\mu''\theta + \frac 6R \mu' \theta' \right] \partial_x \tilde v + \left[\frac{\mu'''}{R^3}\theta + \frac3{R^2}\mu'' \theta' + \frac3R \mu' \theta'' + \mu'\phi'\theta \right] \tilde v,\text{ and}\label{3.10}\\
F_3&:=\gamma \mu \theta[\overline u_1 \partial_x \tilde u + u_2 \partial_x \overline{\tilde u} + \partial_x \overline u_1 \tilde u + \partial_x u_2 \overline{\tilde u}].\label{3.11}
\end{align}

We define $\psi$ and $D$ as in Lemma \ref{le_schr}. Taking into account the supports of the derivatives of $\theta$ and the derivatives of $\mu$, it is easy to see that
\begin{align}
\|e^\psi(-\mu\theta(v_1+v_2)\partial_x \tilde v)\|_{L^2}& \leq C \|e^\psi \partial_x g\|_{L^2(D)} + C e^{9\alpha} \|\tilde v\|_{L^2(D)} + C e^{25\alpha} A_{k,R}(\tilde v),\label{3.12}\\
\|e^\psi (-\mu\theta \partial_x(v_1+v_2)\tilde v)\|_{L^2(D)} &\leq C \|e^\psi g\|_{L^2(D)},\label{3.13}\\
\|e^\psi F_1\|_{L^2(D)}&\leq Ce^{25\alpha} A_{k,R} (\tilde v), \label{3.14}\\
\|e^\psi F_2\|_{L^2(D)}&\leq C e^{9\alpha} (\|\tilde v \|_{L^2(D)}+\|\partial_x \tilde v\|_{L^2(D)} + \|\partial_x^2 \tilde v \|_{L^2(D)}),\text{ and}\label{3.15}\\
\notag\|e^\psi F_3\|_{L^2(D)}&\leq \widetilde C \|e^\psi \mu \theta \partial_x \tilde u\|_{L^2(D)} + \widetilde C \|e^\psi \mu \theta \tilde u\|_{L^2(D)}\\
&\leq \widetilde C \|e^\psi \partial_x \tilde g\|_{L^2(D)} + \widetilde C e^{9\alpha} \|\tilde u\|_{L^2(D)} + \widetilde C e^{25\alpha} A_{s,R}(\tilde u) + \widetilde C \| e^\psi \tilde g\|_{L^2(D)},\label{3.16}
\end{align}

where the constant $\widetilde C$ is independent of $\|\phi'\|_{L^\infty}$.\\

Applying Lemma \ref{le_kdv} to the function $g$ and taking into account \eqref{3.8} to \eqref{3.16} we have that for $\overline C:=\max\{\|\phi'\|^2_{L^\infty},\|\phi''\|_{L^\infty}\}=\overline C(r)>0$ and $\alpha\geq \overline C R^{3/2}$
\begin{align}
\notag \frac{\alpha^{5/2}}{R^3} \|e^\psi g\|_{L^2(D)} + \frac{\alpha^{3/2}}{R^2} \|e^\psi \partial_x g\|_{L^2(D)} \leq &C \|e^\psi g\|_{L^2(D)} + C \|e^\psi \partial_x g\|_{L^2(D)} + C e^{25\alpha} A_{k,R}(\tilde v) + C e^{9\alpha}\\
& + \widetilde C \|e^\psi \tilde g\|_{L^2(D)} + \widetilde C \|e^\psi \partial_x \tilde g\|_{L^2(D)} + \widetilde C e^{25\alpha} A_{s,R}(\tilde u),\label{3.17}
\end{align}

where constants $C$ and $\widetilde C$ do not depend upon $R$ and $\widetilde C$ is independent of $\|\phi'\|_{L^\infty}$.\\

Now we apply Lemma \ref{le_schr} to the function $\tilde g$. From the definition of $\tilde g$ and the fact that $\tilde u$ satisfies the equation
\begin{align}
i\tilde u_t + \partial_x^2 \tilde u = \beta v_1 \tilde u + \beta u_2 \tilde v - |u_1|^2 \tilde u - u_2 \overline u_1 \tilde u - u_2^2 \overline{\tilde u},\label{3.17b}
\end{align}

it can be seen that
\begin{align}
i\partial_t \tilde g + \partial_x^2 \tilde g = -\mu \theta |u_1|^2 \tilde u- \mu\theta u_2 \overline u_1 \tilde u - \mu \theta u_2^2 \overline{\tilde u} + G_1 + G_2 +G_3, \label{3.18}
\end{align}

where
\begin{align}
G_1&:= \mu [\theta''\tilde u + \theta' \partial_x\tilde u], \label{3.19}\\
G_2&:=2\frac{\mu'}R \theta \partial_x \tilde u + \left[\frac1{R^2} \mu'' \theta +\frac2R\mu'\theta' + i \mu' \phi' \theta \right] \tilde u,\text{ and}\label{3.20}\\
G_3&:=\beta \mu \theta [v_1 \tilde u + u_2 \tilde v]. \label{3.21}
\end{align}

Taking again into account the supports of the derivatives of $\theta$ and $\mu$, it is easy to see that
\begin{align}
\|e^\psi (-\mu \theta |u_1|^2 \tilde u - \mu\theta u_2 \overline u_1 \tilde u - \mu \theta u^2_2 \overline{\tilde u})\|_{L^2(D)} &\leq C \| e^\psi \tilde g \|_{L^2(D)},\label{3.22}\\
\|e^\psi G_1\|_{L^2(D)}&\leq C e^{25\alpha} A_{s,R}(\tilde u),\label{3.23}\\
\|e^\psi G_2\|_{L^2(D)}&\leq C e^{9\alpha}(\|\tilde u\|_{L^2(D)}+\|\partial_x \tilde u\|_{L^2(D)}),\text{ and}\label{3.24}\\
\|e^\psi G_3\|_{L^2(D)}&\leq C \|e^\psi \tilde g\|_{L^2(D)} + C \|e^\psi g\|_{L^2(D)}.\label{3.25}
\end{align}

Applying Lemma \ref{le_schr} to the function $\tilde g$ and taking into account \eqref{3.18} to \eqref{3.25} we have that for $\overline C:=\max\{\|\phi'\|^2_{L^\infty}, \|\phi''\|_{L^\infty}\}=\overline C(r)$ and $\alpha\geq \overline C R^2$
\begin{align}
\frac{\alpha^{3/2}}{R^2} \|e^\psi \tilde g\|_{L^2(D)} + \frac{\alpha^{1/2}}{R} \|e^\psi \partial_x \tilde g\|_{L^2(D)} \leq & C \|e^\psi \tilde g\|_{L^2(D)} + C e^{25\alpha} A_{s,R}(\tilde u) + C e^{9\alpha} + C \|e^\psi g\|_{L^2(D)},\label{3.26}
\end{align}

where the constant $C$ does not depend upon $R$.\\

If we take $\alpha:= \overline C R^{2}$ in \eqref{3.17} and \eqref{3.26}, these two inequalities become
\begin{align}
\notag\overline C^{5/2} R^{2} \|e^\psi  g\|_{L^2(D)} + \overline C^{3/2} R& \|e^\psi \partial_x  g\|_{L^2(D)} \leq  C (\|e^\psi  g\|_{L^2(D)} + \|e^\psi \partial_x g\|_{L^2(D)}+  e^{25\overline C R^{2}} A_{k,R}(\tilde v)\\
& + e^{9\overline C R^{2}}) + \widetilde C (\|e^\psi \tilde g\|_{L^2(D)} + \|e^\psi \partial_x \tilde g\|_{L^2(D)}+ e^{25 C R^{2}} A_{s,R}(\tilde u)),\label{3.27}
\end{align}

and 
\begin{align}
\overline C^{3/2} R \|e^\psi  \tilde g\|_{L^2(D)} + \overline C^{1/2} \|e^\psi \partial_x  \tilde g\|_{L^2(D)} \leq  C (\|e^\psi \tilde g\|_{L^2(D)} +  e^{25\overline C R^{2}} A_{s,R}(\tilde u) + e^{9\overline C R^{2}}+ \|e^\psi  g\|_{L^2(D)} ),\label{3.28}
\end{align}

Adding inequalities \eqref{3.27} and \eqref{3.28}, noting that $\overline C^{1/2}\geq \|\phi'\|_{L^\infty}\geq\frac 1r$, taking $r_1$ small enough in order to have $\frac12 \overline C^{1/2}\geq \widetilde C$ for $r\in(0,r_1)$, the term $\widetilde C \|e^\psi \partial_x \tilde g\|_{L^2(D)}$ on the right hand side of the new inequality can be absorbed by $\frac12 \overline C^{1/2} \|e^{\psi} \partial_x \tilde g\|_{L^2(D)}$ on the left hand side. Also observing that, for sufficiently large $R$, the terms $C\|e^\psi g\|_{L^2(D)}$, $C\|e^\psi \tilde g\|_{L^2(D)}$, and $C\|e^\psi \partial_x g\|_{L^2(D)}$ can be absorbed, respectively, by $\frac12 \overline C^{5/2} R^2 \|e^\psi g\|_{L^2(D)}$, $\frac12 \overline C^{3/2} R\|e^\psi \tilde g\|_{L^2(D)}$, and $\frac12 \overline C^{3/2} R \|e^\psi \partial_x g\|_{L^2(D)}$, we can conclude that\\

\begin{align}
\notag \frac12 \overline C^{5/2} R^{2} \|e^{\psi} g \|_{L^2(D)} +\frac12 \overline C^{3/2} \|e^\psi \partial_x g\|_{L^2(D)} + \frac12 \overline C^{3/2} R \|e^\psi \tilde g\|_{L^2(D)} +\frac12 \overline C^{1/2}\|e^\psi \partial_x \tilde g\|_{L^2(D)}\\
\leq  C (e^{25 \overline C R^{2}} A_R(\tilde u,\tilde v) + e^{9\overline C R^{2}}).\label{3.29}
\end{align}

Since $\tilde g=\tilde u$, and $g=\tilde v$ in $Q_r$, and $\psi(x,t)\geq 16\alpha$ for $(x,t)\in Q_r$, we have that
\begin{align}
\frac12 \overline C^{5/2} R^{2} e^{16 \overline C R^{2}} \|\tilde v \|_{L^2(Q_r)} + \frac12 \overline C^{3/2} R e^{16\overline C R^{2}} \|\tilde u \|_{L^2(Q_r)} \leq C (e^{25\overline C  R^{2}} A_R(\tilde u, \tilde v) + e^{9\overline C R^{2}}),\label{3.30}
\end{align}
for $r\in(0,r_1)$.\\

Taking into account that $\|\tilde v\|_{L^2(Q_r)} + \|\tilde u\|_{L^2(Q_r)}>0$, the term $C e^{9\overline C R^{2}}$, on the right hand side of \eqref{3.30}, can be absorbed by $\frac14 \overline C^{3/2} R e^{16\overline C R^{2}}(\|\tilde v\|_{L^2(Q_r)} + \|\tilde u\|_{L^2(Q_r)})$, when we take $R$ sufficiently large.\\

In this manner we can conclude that there exists $R_0\geq 2$ such that for $R\geq R_0$,
$$\|\tilde u\|_{L^2(Q_r)} + \|\tilde v\|_{L^2(Q_r)}\leq C e^{9\overline C R^{2}} A_R(\tilde u,\tilde v),$$
and Theorem \ref{le_th} is proved.\qed\\

\section{Proof of Theorem \ref{main_th}}\label{SectionMR} In the proof of Theorem \ref{main_th} we will use the following interpolation Lemma (see \cite{NP2009}).

\begin{lemma}\label{interpolation} Let $a,b>0$. Assume that $J^a f:=(1-\partial_x^2)^{a/2} f\in L^2(\mathbb R)$ and $(1+x^2)^{b/2}f\in L^2(\mathbb R)$. Then for any $\theta \in(0,1)$
\begin{align}
\|J^{\theta a}((1+x^2)^{(1-\theta)b/2}f)\|_{L^2(\mathbb R)} \leq C \|(1+x^2)^{b/2} f\|^{1-\theta}_{L^2(\mathbb R)} \|J^a f\|^\theta_{L^2(\mathbb R)}.\label{4.1}
\end{align}

\end{lemma}

\textit{Proof of Theorem \ref{main_th}}. This proof follows the scheme in the proof of Theorem 1.1 in \cite{BIM2013}.\\

Let $\tilde u:=u_1-u_2$ and $\tilde v:=v_1-v_2$. Let $r_1$ be as in Theorem \ref{le_th} and let us take $r$ in $(0,r_1)$. We look for a constant $a_0$, such that if $a>a_0$, then $\tilde u\equiv 0$ and $\tilde v\equiv 0$ in the set $D_0:=\{(x,t):x\in\mathbb R, t\in[r,1-r]\}$. Reasoning by contradiction, suppose that $(\tilde u,\tilde v)\not\equiv (0,0)$ in $D_0$. Without loss of generality we can afirm that if $Q_r:=\{(x,t):0\leq x\leq 1,t\in[ r,1-r]\}$, then
\begin{align}
\|\tilde u\|_{L^2(Q_r)}+\|\tilde v\|_{L^2(Q_r)}>0.\label{4.3}
\end{align}

By Theorem \ref{le_th}, there exist positive constants $\overline C=\overline C(r)$, $C$, and $R_0\geq 2$ such that
\begin{align}
\|\tilde u\|_{L^2(Q_r)} + \|\tilde v\|_{L^2(Q_r)}\leq C e^{9\overline C R^{2}} A_R (\tilde u,\tilde v)\quad \text{for every $R\geq R_0$}. \label{4.4}
\end{align}

For $R>3$, we take $N\in\mathbb N$, with $N>3R$. We construct a $C^\infty$ function $\Phi_{R,N}:\mathbb R\to\mathbb R$ supported in $(R,N+1)$, such that $0\leq \Phi_{R,N}\leq 1$, $\Phi_{R,N}\equiv 1$ in $[R+1,N]$ and $|\Phi_{R,N}^{(j)}(x)|\leq C_j$ for all $x\in\mathbb R$, with $C_j$ independent of $R$ and $N$.\\

Let us define $w$ and $z$ for $(x,t)\in\mathbb R\times [0,1]$ by
\begin{align}
w(t)(x)&:= \Phi_{R,N}(x) \tilde u(t)(x), \text{ and}\label{4.5}\\
z(t)(x)&:= \Phi_{R,N}(x) \tilde v(t)(x).\label{4.6}
\end{align}

Then $w$ and $z$ satisfy the hypotheses of Lemmas \ref{carleman_schr}, and \ref{carleman_kdv}, respectively. Therefore for $\lambda>2$,
\begin{align}
\notag \|e^{\lambda x} w\|_{L^\infty_t L^2_x(D)} + \|e^{\lambda x} \partial_x w\|_{L^\infty_x L^2_t(D)}  \leq & C [\lambda^2 (\|J(e^{\lambda x} w(0))\|_{L^2_x(\mathbb R)} + \|J(e^{\lambda x} w(1))\|_{L^2_x(\mathbb R)})\\
&+\|e^{\lambda x} (iw_t + \partial_x^2 w)\|_{L^1_t L^2_x(D)\cap L^1_x L^2_t(D)}],\label{4.7}\\
\notag\|e^{\lambda x}z\|_{L^\infty_t L^2_x(D)} + \sum_{k=1}^2 \|e^{\lambda x} \partial_x^k z\|_{L^\infty_x L^2_t(D)} \leq & C[\lambda^2 (\| J(e^{\lambda x} z(0))\|_{L^2_x(\mathbb R)} + \|J(e^{\lambda x} z(1))\|_{L^2_x(\mathbb R)})\\
& + \|e^{\lambda x}(z_t+\partial_x^3 z)\|_{L^1_t L^2_x (D) \cap L^1_x L^2_t (D)}].\label{4.8}
\end{align}

But
\begin{align}
i w_t + \partial_x^2 w =& \beta v_1 w + \beta u_2 z - |u_1|^2 w - u_2 \overline u_1 w - u_2^2 \overline w + F,\text{ and}\label{4.9}\\
\notag z_t + \partial_x^3 z=& -\partial_x(v_1+v_2)z - (v_1+v_2)\Phi_{R,N} \partial_x \tilde v\\
&+ \gamma \overline u_1 \Phi_{R,N} \partial_x \tilde u + \gamma u_2 \Phi_{R,N} \partial_x \overline{\tilde u} + \gamma \partial_x u_1 w + \gamma \partial_x u_2 \overline w + G,\label{4.10}
\end{align}

where
\begin{align}
F:=& \Phi''_{R,N} \tilde u + 2 \Phi'_{R,N}\partial_x \tilde u,\text{ and}\label{4.11}\\
G:=&\Phi'''_{R,N} \tilde v + 3 \Phi''_{R,N}\partial_x \tilde v + 3\Phi'_{R,N} \partial^2_x \tilde v.\label{4.11b}
\end{align}

Therefore, from \eqref{4.7} to \eqref{4.10}, and using the facts that $\|\cdot\|_{L^2(D)}\leq \|\cdot\|_{L^\infty_t L^2_x(D)}$ and $\|\cdot\|_{L^1_t L^2_x(D)}\leq C\|\cdot\|_{L^2(D)}$, it follows that
\begin{align}
\|e^{\lambda x}w\|_{L^2(D)} + \|e^{\lambda x} \partial_x w\|_{L^\infty_x L^2_t(D)}\leq I+II+III,\label{4.13}
\end{align}
where
\begin{align*}
I:=& C(\|e^{\lambda x}v_1 w\|_{L^2(D)\cap L^1_x L^2_t(D)} + \|e^{\lambda x} u_2 z\|_{L^2(D)\cap L^1_x L^2_t(D)}+ \|e^{\lambda x}|u_1|^2 w\|_{L^2(D)\cap L^1_x L^2_t(D)}\\
&+ \|e^{\lambda x} u_2 \overline u_1 w \|_{L^2(D)\cap L^1_x L^2_t (D)} + \|e^{\lambda x} u_2^2 \tilde w\|_{L^2(D)\cap L^1_x L^2_t(D)}),\\
II:=& C \|e^{\lambda x} F\|_{L^2(D)\cap L^1_x L^2_t(D)}.\\
III:=& C \lambda^2 (\|J(e^{\lambda x} w(0))\|_{L^2_x(\mathbb R)} + \|J(e^{\lambda x} w(1))\|_{L^2_x(\mathbb R)}),
\end{align*}

and
\begin{align}
\|e^{\lambda x}z\|_{L^2(D)} + \sum_{k=1}^2 \|e^{\lambda x} \partial_x^k z\|_{L^\infty_x L^2_t(D)} \leq \widetilde I + \widetilde {II} + \widetilde {III} + \widetilde {IV}, \label{4.14}
\end{align}

where
\begin{align*}
\widetilde I :=&C(\|e^{\lambda x} \partial_x(v_1+v_2)z\|_{L^2(D)\cap L^1_x L^2_t(D)} + \|e^{\lambda x} \partial_x u_1 w\|_{L^2(D)\cap L^1_x L^2_t (D)} + \|e^{\lambda x}\partial_x u_2 \overline w\|_{L^2(D)\cap L^1_x L^2_t(D)}), \\
\widetilde {II} :=&C (\|e^{\lambda x}(v_1+v_2) \Phi_{R,N} \partial_x \tilde v\|_{L^2(D)\cap L^1_x L^2_t(D)} + \|e^{\lambda x}\overline u_1 \Phi_{R,N}\partial_x \tilde u\|_{L^2(D)\cap L^1_x L^2_t(D)}\\
& + \|e^{\lambda x} u_2 \Phi_{R,N} \partial_x \overline{\tilde u}\|_{L^2(D)\cap L^1_x L^2_t(D)}), \\
\widetilde {III} :=& C \|e^{\lambda x}G\|_{L^2(D)\cap L^1_x L^2_t(D)}, \\
\widetilde {IV} :=& C\lambda^2 (\|J(e^{\lambda x}z(0))\|_{L^2_x(\mathbb R)} + \|J(e^{\lambda x} z(1))\|_{L^2(\mathbb R)} ).
\end{align*}

We now estimate $I,II,III,\widetilde I,\widetilde{II},\widetilde {III}$, and $\widetilde{IV}$ separately.\\

\textit{Estimation of $I$}: By observing that $w$ and $z$ are supported in the set
\begin{align}
D_R:=\{x:x\geq R\}\times [0,1],\label{4.15a}
\end{align}

applying Hölder's inequality, and denoting by $\chi_{D_R}$ the characteristic function of the set $D_R$, we have that
\begin{align}
\notag I\leq & C (\|v_1 \chi_{D_R}\|_{L^\infty(D)} \|e^{\lambda x} w\|_{L^2(D)} + \|v_1 \chi_{D_R}\|_{L^2_x L^\infty_t(D)} \|e^{\lambda x} w\|_{L^2(D)}\\
\notag& + \|u_2 \chi_{D_R}\|_{L^\infty(D)} \|e^{\lambda x}z\|_{L^2(D)} + \|u_2 \chi_{D_R}\|_{L^2_x L^\infty_t (D)} \|e^{\lambda x}z\|_{L^2(D)}\\
\notag& + \| |u_1|^2\chi_{D_R} \|_{L^\infty(D)} \|e^{\lambda x}w\|_{L^2(D)} + \||u_1|^2 \chi_{D_R}\|_{L^2_x L^\infty_t(D)} \|e^{\lambda x}w\|_{L^2(D)}\\
\notag &+ \|u_2 \overline u_1 \chi_{D_R}\|_{L^\infty(D)} \|e^{\lambda x}w\|_{L^2(D)} + \|u_2\overline u_1 \chi_{D_R}\|_{L^2_x L^\infty_t(D)} \|e^{\lambda x} w\|_{L^2(D)}\\
&+\|u_2^2 \chi_{D_R} \|_{L^\infty(D)} \|e^{\lambda x} w\|_{L^2(D)} + \|u_2^2 \chi_{D_R}\|_{L^2_x L^\infty_t (D)} \|e^{\lambda x}w\|_{L^2(D)} ).\label{4.15}
\end{align}

It can be seen that the norms involving $v_1,u_2,|u_1|^2,u_2\overline u_1,$ and $u_2^2$ on the right hand side of \eqref{4.15} tend to zero as $R\to\infty$. Let us illustrate the previous statement in the cases of $\|v_1 \chi_{D_R}\|_{L^\infty(D)}$ and $\|v_1\chi_{D_R}\|_{L^2_x L^\infty_t(D)}$. Using the interpolation Lemma \ref{interpolation}, with $a=3$, $b=3$ and $\theta=\frac15$,
\begin{align*}
\|(1+x^2)^{\frac45\frac32} v_1\|_{H^{3/5}} \leq C \|(1+x^2)^{3/2} v_1\|_{L^2}^{4/5} \|v_1\|^{1/5}_{H^3}.
\end{align*}

Since $H^{3/5}\hookrightarrow L^\infty(\mathbb R)$, then there exists $C>0$, such that
\begin{align*}
\|(1+x^2)^{6/5} v_1\|_{L^\infty(\mathbb R_x)} \leq C \|(1+x^2)^{3/2} v_1\|^{4/5}_{L^2(\mathbb R_x)} \|v_1\|^{1/5}_{H^3(\mathbb R_x)}\leq C.
\end{align*}

In this manner
\begin{align*}
(1+x^2)^{6/5} |v_1(x,t)|&\leq C \text{ for all $x\in\mathbb R$ and $t\in[0,1]$ and},\\
|v_1(x,t) \chi_{D_R} (x,t)|&\leq \frac C{(1+R^2)^{6/5}} \text{ for all $x\in\mathbb R$ and $t\in[0,1]$ and},
\end{align*}

in consequence $\|v_1 \chi_{D_R}\|_{L^\infty(D)}\to 0$ as $R\to\infty$.\\

On the other hand,
\begin{align*}
\|v_1\chi_{D_R} \|_{L^2_x L^\infty_t (D)} = \left( \int_{\mathbb R} \|v_1(x,\cdot_t) \chi_{D_R}(x,\cdot_t)\|_{L^\infty_t[0,1]}^2 dx \right)^{1/2}.
\end{align*}

Since for all $R>0$ and all $x\in \mathbb R$
\begin{align*}
\|v_1(x,\cdot_t) \chi_{D_R} (x,\cdot_t)\|_{L^\infty_t([0,1])}^2 \leq \frac C{(1+x^2)^{12/5}}\in L^1(\mathbb R_x),
\end{align*}

and
\begin{align*}
\lim_{R\to\infty} \|v_1(x,\cdot_t) \chi_{D_R}(x,\cdot_t)\|_{L^\infty_t([0,1])}^2=0,
\end{align*}

then, by the dominated convergence theorem, it follows that $\|v_1\chi_{{D_R}}\|_{L^2_x L^\infty_t(D)}\to 0$ as $R\to\infty$.\\

\textit{Estimation of $\widetilde I$}: As in the estimation of $I$, we have that
\begin{align}
\notag \widetilde I \leq & C ( \|\partial_x (v_1 + v_2)\chi_{D_R}\|_{L^\infty(D)} + \|\partial_x (v_1 + v_2)\chi_{D_R}\|_{L^2_x L^\infty_t(D)}) \|e^{\lambda x} z\|_{L^2(D)}\\
\notag & + C (\|\partial_x u_1 \chi_{D_R}\|_{L^\infty (D)} + \|\partial_x u_1 \chi_{D_R}\|_{L^2_x L_t^\infty(D)}) \|e^{\lambda x} w \|_{L^2(D)} \\
& + C (\|\partial_x u_2 \chi_{D_R}\|_{L^\infty (D)} + \|\partial_x u_2 \chi_{D_R}\|_{L^2_x L_t^\infty(D)}) \| e^{\lambda x}w\|_{L^2(D)}. \label{4.16}
\end{align}

As in the estimation of $I$, it can be also seen, using Lemma \ref{interpolation} and the dominated convergence theorem, that the norms involving $\partial_x(v_1+v_2)$, $\partial_x u_1$, and $\partial_x u_2$ on the right hand side of \eqref{4.16} tend to zero as $R\to \infty$.\\

In this way, adding the inequalities \eqref{4.13} and \eqref{4.14}, the terms $I$ and $\widetilde I$, on the right hand side of the new inequality, can be absorbed, for $R$ large enough, by the terms $\|e^{\lambda x} w\|_{L^2(D)}$ and $\|e^{\lambda x}z\|_{L^2(D)}$ on the left hand side of this inequality, so for $R$ large enough, we have that
\begin{align}
\notag\|e^{\lambda x}w\|_{L^2(D)} + \|e^{\lambda x} \partial_x w\|_{L_x^\infty L_t^2(D)} + \|e^{\lambda x} z\|_{L^2(D)} + \sum_{k=1}^2 \|e^{\lambda x} \partial^k_x z\|_{L^\infty_x L^2_t(D)}\\
\leq II + III + \widetilde{II} + \widetilde{III} + \widetilde{IV}. \label{4.17}
\end{align}

\textit{Estimation of $\widetilde{II}$}: Proceeding as before,
\begin{align}
\notag\widetilde{II}=& C ( \|e^{\lambda x}(v_1+v_2) [\partial_x (\Phi_{R,N} \tilde v)]-\Phi'_{R,N} \tilde v\|_{L^2(D)\cap L^1_xL^2_t(D)}\\
\notag &+ \|e^{\lambda x} \overline u_1[\partial_x (\Phi_{R,N} \tilde u) - \Phi'_{R,N} \tilde u]\|_{L^2(D)\cap L^1_x L^2_t(D)} + \|e^{\lambda x} u_2[\partial_x(\Phi_{R,N} \overline{\tilde u}) - \Phi'_{R,N} \overline{\tilde u}]\|_{L^2(D)\cap L^1_xL^2_t(D)})\\
\notag\leq &C[ (\|(v_1+v_2)\chi_{D_R}\|_{L^2_xL^\infty_t(D)} + \|(v_1+v_2)\chi_{D_R}\|_{L^1_xL^\infty_t(D)}) \|e^{\lambda x} \partial_x z\|_{L^\infty_x L^2_t(D)} \\
\notag &+ (\|(v_1+v_2)\chi_{D_R}\|_{L^\infty (D)} + \|(v_1+v_2)\chi_{D_R}\|_{L^2_xL^\infty_t(D)}) \|e^{\lambda x} \Phi'_{R,N} \tilde v\|_{L^2(D)}\\
\notag & + (\|\overline u_1 \chi_{D_R}\|_{L^2_x L^\infty_t(D)} + \|\overline u_1 \chi_{D_R}\|_{L^1_x L^\infty_t(D)}) \|e^{\lambda x} \partial_x w\|_{L^\infty_x L^2_t(D)}\\
\notag & + (\|\overline u_1 \chi_{D_R}\|_{L^\infty(D)} + \|\overline u_1 \chi_{D_R}\|_{L^2_x L^\infty_t(D)}) \|e^{\lambda x} \Phi'_{R,N} \tilde u\|_{L^2(D)}\\
\notag & + (\|u_2\chi_{D_R}\|_{L^2_x L^\infty_t(D)} + \|u_2 \chi_{D_R}\|_{L^1_x L^2_t(D)}) \|e^{\lambda x} \partial_x \overline w\|_{L^\infty_x L^2_t(D)}\\
&+(\|u_2 \chi_{D_R}\|_{L^\infty(D)} + \|u_2 \chi_{D_R}\|_{L^2_x L^\infty_t(D)}) \|e^{\lambda x} \Phi'_{R,N} \overline {\tilde u}\|_{L^2(D)}].\label{4.18}
\end{align}

The norms involving $(v_1+v_2),\overline u_1,\overline u_2$ on the right hand side of \eqref{4.18} tend to zero as $R\to\infty$. In this way, the terms on the right hand side of \eqref{4.18}, including the factors $\|e^{\lambda x} \partial_x z\|_{L^\infty_x L^2_t(D)}$, $\|e^{\lambda x}\partial_x w\|_{L^\infty_x L^2_t(D)}$ and $\|e^{\lambda x}\partial_x \overline w\|_{L^\infty_x L^2_t(D)}$, can be absorbed, for $R$ large enough, by the terms $\|e^{\lambda x} \partial_x z\|_{L^\infty_x L^2_t (D)}$ and $\|e^{\lambda x} \partial_x w\|_{L^\infty_x L^2_t(D)}$ on the left hand side of \eqref{4.17}. In consequence, from \eqref{4.17}, \eqref{4.18}, and the previous observation, it follows that, for $R$ large enough,
\begin{align}
\notag \|e^{\lambda x} &w\|_{L^2(D)} + \|e^{\lambda x} \partial_x w\|_{L^\infty_x L^2_t(D)} + \|e^{\lambda x} z\|_{L^2(D)} + \sum_{k=1}^2 \|e^{\lambda x} \partial_x^k z\|_{L^\infty_x L^2_t(D)}\\
\notag \leq & C [ (\|(v_1+v_2)\chi_{D_R}\|_{L^\infty(D)} + \|(v_1+v_2)\chi_{D_R}\|_{L^2_x L^\infty_t (D)}) \|e^{\lambda x} \Phi'_{R,N} \tilde v\|_{L^2(D)}\\
\notag &+ (\|\overline u_1 \chi_{D_R}\|_{L^\infty(D)} + \|\overline u_1 \chi_{D_R}\|_{L^2_x L^\infty_t(D)}) \|e^{\lambda x} \Phi'_{R,N}\tilde u\|_{L^2(D)}\\
& (\|u_2 \chi_{D_R}\|_{L^\infty(D)} + \|u_2\chi_{D_R}\|_{L^2_x L^\infty_t(D)}) \|e^{\lambda x} \Phi'_{R,N} \overline{\tilde u}\|_{L^2(D)} ] + II + III + \widetilde{III} + \widetilde{IV}. \label{4.19} 
\end{align}

Taking into account that $\tilde u(0),\tilde u(1),\tilde v(0),\tilde v(1)\in H^1(e^{a x^{2}}dx)$, then for all $\overline \lambda>0$, $\tilde u(0)$, $\tilde u(1)$, $\tilde v(0)$, $\tilde v(1)\in H^1(e^{2\overline \lambda |x|} dx)$. In this manner, applying to the equations \eqref{3.17b}, and \eqref{3.7}, satisfied by $\tilde u$ and $\tilde v$, respectively, a standard procedure, similar to that given in Theorem 1.3 of \cite{BIM2011} we can observe that $\tilde u$ and $\tilde v$ are bounded functions from the time interval $[0,1]$ with values in $H^1(e^{2\overline\lambda |x|}dx)$.\\

If we take $\overline \lambda:=\lambda +1$, it follows that
\begin{align}
\notag \|e^{\lambda x} \Phi'_{R,N}\tilde v\|_{L^2(D)} \leq & C \|e^{\lambda x} \Phi'_{R,N}\tilde v\|_{L^2([R,R+1]\times [0,1])} + C \|e^{\lambda x} \Phi'_{R,N}\tilde v\|_{L^2([N,N+1]\times [0,1])}\\
\notag \leq &C e^{\lambda (R+1)}\|\Phi'_{R,N}\tilde v\|_{L^2(D)} + C \|e^{-x} e^{\overline \lambda x}\Phi'_{R,N}\tilde v\|_{L^2([N,N+1]\times[0,1])}\\
\notag \leq & C e^{\lambda(R+1)} + C e^{-N} \|e^{\overline \lambda |x|} \Phi'_{R,N} \tilde v\|_{L^2([N,N+1]\times [0,1])}\\
\leq &C e^{\lambda (R+1)} + C e^{-N}\|e^{\overline \lambda |x|}\tilde v\|_{L^2(D)}\leq  C e^{\lambda (R+1)} + C_\lambda e^{-N}.\label{4.20}
\end{align}

In a similar manner we also have that
\begin{align}
\|e^{\lambda x} \Phi'_{R,N} \tilde u\|_{L^2(D)},\|e^{\lambda x}\Phi'_{R,N}\overline{\tilde u}\|_{L^2(D)} \leq C e^{\lambda(R+1)}+C_\lambda e^{-N}.\label{4.21}
\end{align}

This way, from \eqref{4.19}, \eqref{4.20} and \eqref{4.21}, it follows that
\begin{align}
\notag \|e^{\lambda x} w\|_{L^2(D)} + \|e^{\lambda x} \partial_x w \|_{L^\infty_x L^2_t(D)} + &\|e^{\lambda x} z\|_{L^2(D)} + \sum_{k=1}^2 \|e^{\lambda x} \partial_x^k z\|_{L^\infty_x L^2_t(D)}\\
\leq & C e^{\lambda(R+1)} + C_\lambda e^{-N} + II + III + \widetilde{III} + \widetilde{IV}.\label{4.22}
\end{align}

\textit{Estimation of $II$ and $\widetilde {III}$}: Reasoning in a similar way to how we did in the estimate of $\widetilde {II}$, it can be seen that
\begin{align}
II & = \|e^{\lambda x}F\|_{L^2(D)\cap L^1_xL^2_t(D)} \leq C e^{\lambda(R+1)} + C_\lambda e^{-N},\text{ and} \label{4.23}\\
\widetilde{III} & = \|e^{\lambda x} G\|_{L^2(D)\cap L^1_x L^2_t(D)} \leq C e^{\lambda(R+1)} + C_\lambda e^{-N}.\label{4.24}
\end{align}

From \eqref{4.22}, \eqref{4.23} and \eqref{4.24} it follows that
\begin{align}
\notag \|e^{\lambda x} w\|_{L^2(D)} + \|e^{\lambda x}\partial_x w\|_{L^\infty_x L^2_t(D)} +& \|e^{\lambda x}z\|_{L^2(D)} + \sum_{k=1}^2 \|e^{\lambda x} \partial_x^k z\|_{L^\infty_x L^2_t(D)}\\
\leq & C e^{\lambda(R+1)} + C_\lambda e^{-N} + III + \widetilde{IV}.\label{4.25}
\end{align}

\textit{Estimation of $III$ and $\widetilde {IV}$}: Since $\Phi_{R,N}$ and its derivatives are bounded by a constant independent of $R$ and $N$, it follows that
\begin{align}
\notag III \leq & C \lambda^2(\lambda+1)(\|e^{\lambda x} \tilde u(0)\|_{L^2(\{x\geq R\})} + \|e^{\lambda x} \partial_x \tilde u(0)\|_{L^2(\{x\geq R\})}\\
& + \|e^{\lambda x} \tilde u(1)\|_{L^2(\{x\geq R\})} + \|e^{\lambda x} \partial_x \tilde u(1)\|_{L^2(\{x\geq R\})})=: H_R(\lambda),\text{ and} \label{4.26}\\
\notag \widetilde {IV} \leq & C \lambda^2(\lambda+1)(\|e^{\lambda x}\tilde v(0)\|_{L^2(\{x\geq R\})}+ \|e^{\lambda x}\partial_x \tilde v(0)\|_{L^2(\{x\geq R\})}\\
& + \|e^{\lambda x}\tilde v(1)\|_{L^2(\{x\geq R\})} + \|e^{\lambda x}\partial_x \tilde v(1)\|_{L^2(\{x\geq R\})})=:\widetilde H_R(\lambda).\label{4.27}
\end{align}

We now consider the region $\Omega_R:=\{(x,t): x\in[3R-1,3R],t\in[0,1]\}$. Then $\Omega_R\subset D_R$, and since $N>3R$, $\tilde u$ and $w$ coincide in $\Omega_R$, and $\tilde v$ and $z$ coincide in $\Omega_R$. We now return to \eqref{4.25}, replace its left hand side by a smaller amount, apply \eqref{4.26} and \eqref{4.27} and make $N\to\infty$ to obtain that
\begin{align}
\notag\|e^{\lambda x}\tilde u\|_{L^2(\Omega_R)} + \|e^{\lambda x}\partial_x \tilde u\|_{L^2(\Omega_R)} + & \|e^{\lambda x}\tilde v\|_{L^2(\Omega_R)} + \sum_{k=1}^2 \|e^{\lambda x} \partial_x^k\tilde v\|_{L^2(\Omega_R)}\\
\leq & C e^{\lambda(R+1)} + H_R(\lambda) + \widetilde H_R(\lambda).\label{4.28}
\end{align}

For $(x,t)\in\Omega_R$,
\begin{align}
\lambda x\geq \lambda(3R-1)\geq 2\lambda R.\label{4.29}
\end{align}

In this way, bearing in mind the definition of $A_{3R}(\tilde u,\tilde v)$, given in the statement of Theorem \ref{le_th}, we conclude from \eqref{4.28} and \eqref{4.29} that
\begin{align}
e^{2\lambda R} A_{3R}(\tilde u,\tilde v)\leq C e^{\lambda(R+1)} + H_R(\lambda) + \widetilde H_R(\lambda).\label{4.30}
\end{align}

For $x\geq R$ and $\lambda:=\frac1{10}aR$, we have that $\lambda x\leq \frac1{10}ax^{2}$, and in consequence, for $R$ large enough,
\begin{align}
\notag H_R\left(\frac1{10}aR\right) \leq & C a^3 R^3(\|e^{\frac1{10}ax^{2}} \tilde u(0)\|_{L^2(\{x\geq R\})} + \|e^{\frac1{10}ax^{2}} \partial_x \tilde u(0)\|_{L^2(\{x\geq R\})} \\
\notag &+ \|e^{\frac1{10}ax^{2}} \tilde u(1)\|_{L^2(\{x\geq R\})} + \|e^{\frac1{10}ax^{2}} \partial_x \tilde u(1)\|_{L^2(\{x\geq R\})})\\
\notag \leq & C a^3 (\|e^{\frac1{2}ax^{2}} \tilde u(0)\|_{L^2(\{x\geq R\})} + \|e^{\frac1{2}ax^{2}} \partial_x \tilde u(0)\|_{L^2(\{x\geq R\})} \\
\notag &+ \|e^{\frac1{2}ax^{2}} \tilde u(1)\|_{L^2(\{x\geq R\})} + \|e^{\frac1{2}ax^{2}} \partial_x \tilde u(1)\|_{L^2(\{x\geq R\})})\\
\leq & C a^3 (\|\tilde u(0)\|_{H^1(e^{ax^{2}}dx)} + \|\tilde u(1)\|_{H^1(e^{ax^{2}}dx)}) \equiv C_a,\label{4.31}
\end{align}
and, similarly,
\begin{align}
\widetilde H_R\left(\frac1{10}aR\right)\leq \widetilde C_a.\label{4.32}
\end{align}

From \eqref{4.30}, \eqref{4.31}, and \eqref{4.32}, taking $\lambda:=\frac1{10}a R$, it follows that for $R$ large enough,
\begin{align*}
e^{\frac15aR^{2}}A_{3R} (\tilde u,\tilde v)\leq C e^{\frac2{15}a R^{2}} + \overline C_a \leq C_a e^{\frac2{15} aR^{2}},
\end{align*}

and thus
\begin{align}
A_{3R} (\tilde u,\tilde v)\leq C_a e^{-\frac 1{15}a R^{2}}.\label{4.33}
\end{align}

From \eqref{4.4} and \eqref{4.33}, it follows that for $R$ large enough
\begin{align*}
\|\tilde u\|_{L^2(Q_r)} + \|\tilde v\|_{L^2(Q_r)} \leq C e^{9\overline C(3R)^{2}} A_{3R}(\tilde u,\tilde v)\leq C_a e^{9\overline C (3R)^{2}} e^{-\frac1{15} a R^{2}}.
\end{align*}

Hence, if $\frac1{15} a > 81\overline C $, i.e., $a>a_0:=1215 \overline C $, then, making $R$ tend to infinity, we conclude that $\|\tilde u\|_{L^2(Q_r)}+\|\tilde v\|_{L^2(Q_r)}=0$, which contradicts the original fact that $\|\tilde u\|_{L^2(Q_r)} + \|\tilde v\|_{L^2(Q_r)}>0$. Therefore $(\tilde u,\tilde v)=(0,0)$ in $D_0$, and Theorem \ref{main_th} is proved.\qed\\

\textbf{Acknowledgments}\\

This work was partially supported by Universidad Nacional de Colombia, Sede-Medellín– Facultad de Ciencias – Departamento de Matemáticas – Grupo de investigación en Matemáticas de la Universidad Nacional de Colombia Sede Medellín, carrera 65 No. 59A -110, post 50034, Medellín Colombia. Proyecto: Análisis no lineal aplicado a problemas mixtos en ecuaciones diferenciales parciales, código Hermes 60827. Fondo de Investigación de la Facultad de Ciencias empresa 3062.

\end{document}